\documentclass[11pt]{article} 
\usepackage{hfstyle}
\usepackage{graphicx}
\usepackage[T1]{fontenc}
\usepackage{hfstyle}
 
\begin{document}
\title{The Tale of Two Queens and Two Towering Figures}
\author{Henryk Fuk\'s
      \oneaddress{
         Department of Mathematics and Statistics, Brock University\\
         St. Catharines, Ontario  L2S 3A1, Canada\\
         Toronto, Ontario M5T 3J1, Canada\\
         \email{hfuks@brocku.ca}
       }
   }

%
\date{}
\Abstract{
In the first article in the series examining mathematics on coins, we discuss
two great scientists who are not only featured on many coins, but also contributed
to both theory and practice of coinage. The first one is Nicolaus Copernicus, author of 
the treaty \emph{Monetae cudendae ratio} in which he proposed the law governing competition between money
and proposed reform of coinage in the Royal Prussia. The second one
is Sir Isaac Newton, warden and Master of the Royal Mint, and overseer of the
Great Recoinage of 1696.
}
\maketitle

\section{Two queens}
Almost two thousand years ago, Roman writer Gaius Petronius Arbiter (ca. 27--66 AD) in his 
novel ``Satiricon'' made a cynical remark about money,
\emph{sola pecunia regnat} -- only money rules (Satyricon XIV). Although  this phrase is taken out of its proper context, 
it became widely known  in the form of a frequently quoted  
proverb, \emph{pecunia regina mundi} -- ``money is the queen of the world''.
Eighteen centuries later, Carl Friedrich Gauss (1777--1855), as his biographer W. Sartorius von Waltershausen 
wants us to believe, uttered words which became one of the most famous mathematical quotes,
\emph{die Mathematik ist die K\"onigin der Wissenschaften} --  ``mathematics is the queen of sciences''.

%
Thus we have two queens. Their domains sometimes overlap, and this overlap will be the subject of this article.
There are at least two levels on which you can consider connections between mathematics and money.
In the second half of the 20-th century, many commemorative coins featuring famous mathematicians have been issued
by many countries in Europe and elsewhere. Collecting and studying them can be a great joy, not to mention that
one can learn along the way a lot of interesting things about the history of the queen of sciences.

There is, however, yet another connection. Some mathematicians took active interest in minting of coins
and in the monetary policy, making both theoretical and practical contributions to these fields. It happens that two of the greatest
mathematical scientists of all times,  Nicolaus Copernicus and Isaac Newton, were both entrusted with duties related to
coinage reforms. In what follows, we will describe ``what they did for coins and what coins did for them''.
%
%
%
%
%
%
%
%

Before we proceed, we need to make some clarification as to whom we consider to be a mathematician.
One has to understand that in the past the distinction between various fields of intellectual endeavor
was not as sharp as it is today. People who made important contributions to mathematics
were often primarily active in other disciplines utilizing mathematics, especially in fields which we call 
today  physics and astronomy. For the purpose of this article, therefore, we will call them all 
\emph{mathematical scientists}.

\section{Nicolaus Copernicus}
Nicolaus Copernicus (1473 -- 1543) is best known for his work \emph{De revolutionibus orbium coelestium}
and creation of the heliocentric system. He interests, however, went far beyond astronomy alone, including
such diverse fields as law, medicine, diplomacy and politics, as well as economics. He 
was a Catholic cleric and a canon at Frombork Cathedral who actively participated in economic administration
of properties belonging to the cathedral chapter, and even directed the defense of Olsztyn
against Teutonic Knights during the Polish--Teutonic War in 1521.

Many coins commemorating his life and achievements have been produced 
in the past, and in order to describe them in chronological order, we have to go back to the years 
following World War~I.

When Poland regained independence in 1918 after 123 years of partitions, one of the gravest and most pressing	 problems it
faced was the monetary chaos and galloping inflation. In order to solve this problem, a new monetary system has been introduced
in 1924, based on the gold standard. In the spirit of continuity with the former Polish-Lithuanian Commonwealth,
the unit of the new currency was to be called z\l{}oty (Eng. ``golden''), equivalent
to $9/31$ gram of pure gold.  Newly reestablished Polish Mint in Warsaw experimented with many different designs
for all coins  prescribed by the monetary reform bill, including the coin with the highest denomination of $100$ z\l{}. 
In 1925, a pattern coin for $100$ z\l{} denomination featuring Copernicus had been produced.
Its reverse  depicted portrait of the astronomer holding some sort of  astronomical
instrument, perhaps a part of armillary sphere or astronomical rings. Several versions of this pattern
had been issued, some carrying edge inscription  \emph{salus republicae suprema lex} - ``prosperity of the republic 
is the supreme law``, others with smooth edge. All are extremely rare now and fetch very high prices when they appear on auctions. 
\begin{center}
  \includegraphics[width=4.5cm]{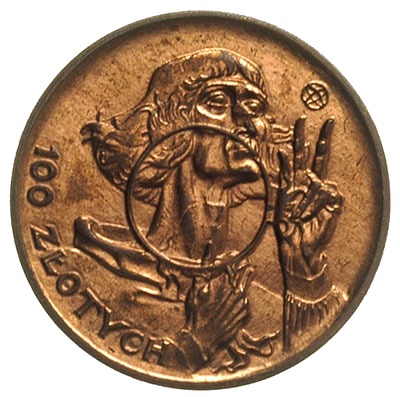} \includegraphics[width=4.5cm]{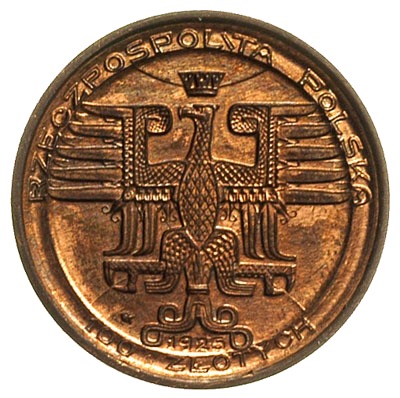}
\end{center}

Unfortunately, coins with denomination 100 z\l{} (as well as 50 z\l{}) have never entered circulation. Shortages of precious metals,
limited capacity of the mint, and huge demand for smaller denominations had been the main cause.
When the war catastrophe of 1939 has prematurely ended the emission of Polish z\l{}oty, the highest-denominated coin
remaining in circulation was 20 z\l{} gold coin with the effigy of Boles\l{}aw I Chrobry, the first King of Poland. The Warsaw mint 
had been plundered by the invading German army, and in 1944, shortly after the eruption of the Warsaw Uprising, retreating Germans
destoyed the mint building. Copernicus had to wait a few more decades before his portrait had a chance to appear on a circulating
coin.

This happened in 1959, when the first coin featuring Copernicus  entered the circulation in the  ``liberated'' Poland. 
It carried the face value of $10$  z\l{},
but the new inscription, People's Republic of Poland, made it clear that the currency, established by the
communist government, had nothing to do with the pre-war z\l{}oty based on gold standard. Nevertheless, the coin, designed by 
an accomplished sculptor J\'ozef Gos\l{}awski and made
in cupronickel, is, in my opinion, the nicest Copernicus coin which has ever been produced. Its only problem was rather big
size (31 mm) and weight (12.9 g), thus in 1967 the mint reduced it slightly, to 28 mm of diameter and 9.5 g of weight. The total
emission of both the big and small version of this coin in years 1959-1969 exceeded 35 mln, and \emph{dycha z Kopernikiem}
(affectionate Polish term for ``ten with Copernicus'') remained in circulation well into 80's. When I was growing up,
one needed two such coins to buy a bar of chocolate, thus in my mind it will for ever be associated with chocolate.
\begin{center}
  \includegraphics[width=4.5cm]{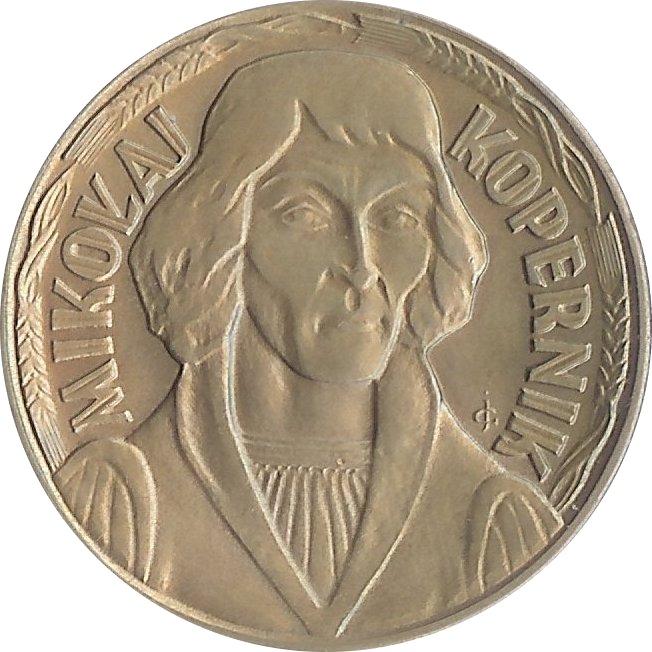} \includegraphics[width=4.5cm]{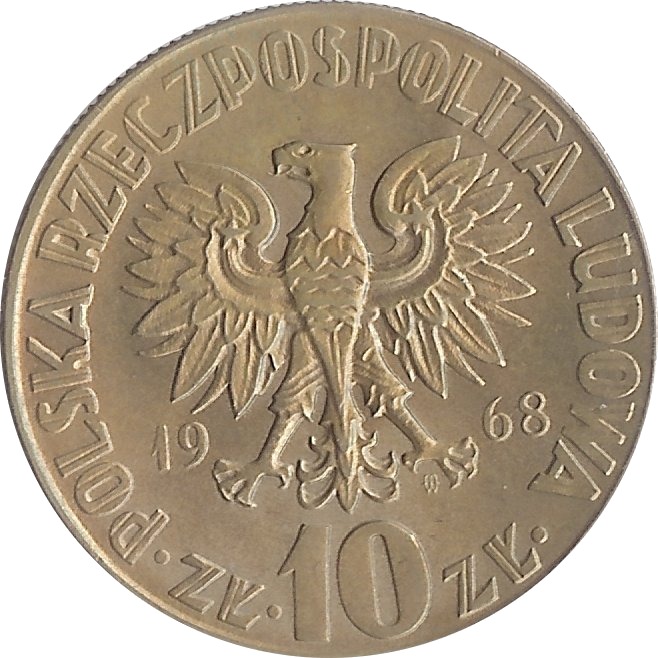}
 \end{center}

In 1973, a big Copernican anniversary was celebrated, the 500-th anniversary of the birth of Copernicus. Both
Poland and Germany issued commemorative coins on this occasion. The Polish coin, carrying 100 z\l{} face value, was
made of 0.625 silver. It was intended for collectors only, and 101 000 pieces have been minted in 1973 and 1974.
The German coin, also made of 0.625 silver, was smaller and lighter, but it was intended for circulation.
8 mln pieces of this coin had been produced, as a part of a larger series of 5 DM silver commemorative coins issued in years
1952--1986.
\begin{center}
 \includegraphics[width=4.0cm]{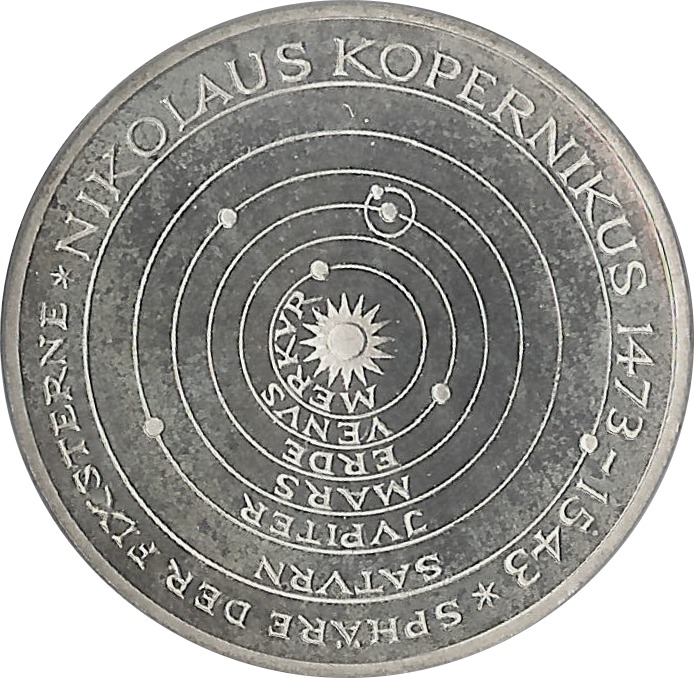} \includegraphics[width=4.0cm]{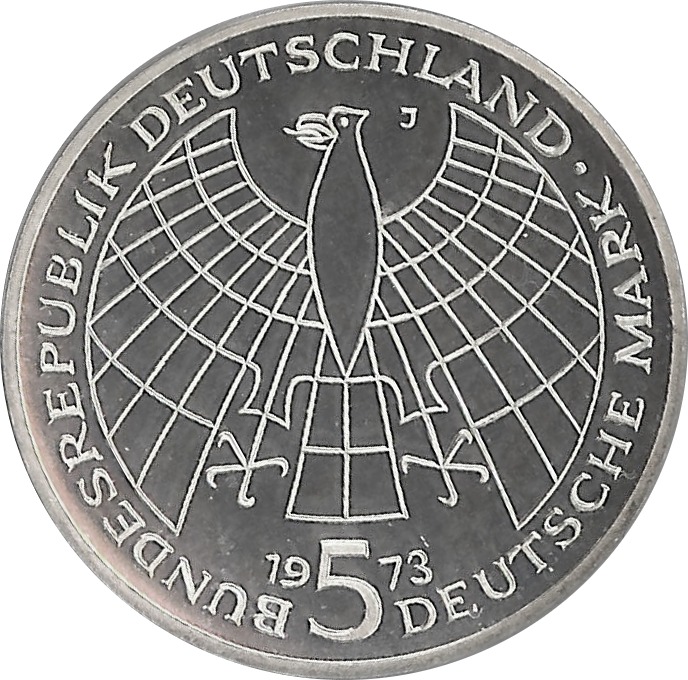}
\includegraphics[width=4.0cm]{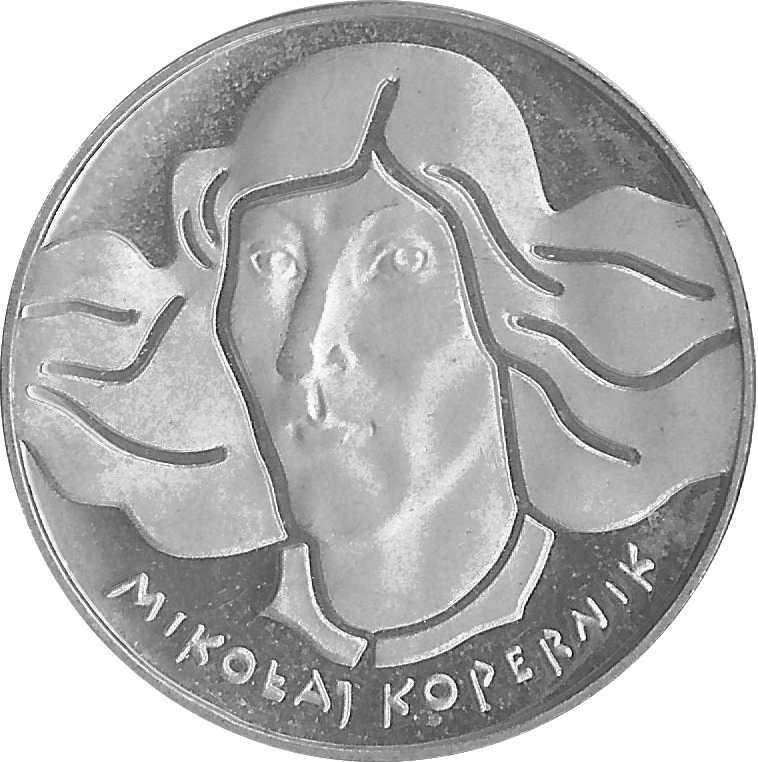} \includegraphics[width=4.0cm]{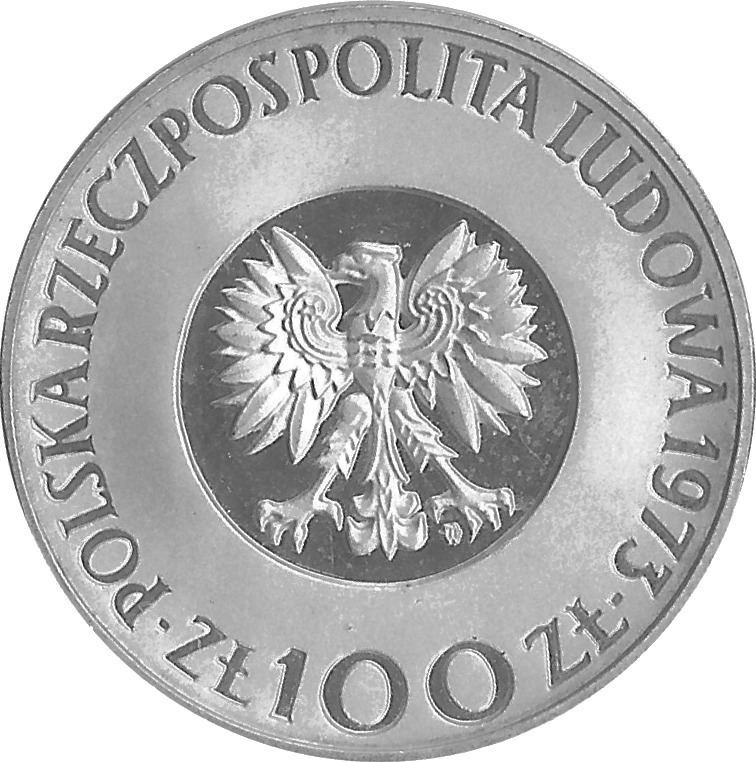}
\end{center}
Its reverse shows the diagram of the solar system as it appeared in \emph{De revolutionibus}, except that the
legends are in German. The diagram, as in the original, shows six concentric circles, representing orbits
of planets from Mercury to Saturn. The seventh and the largest circle is the edge of the coin, here labelled \emph{Sph\"are der Fixsterne}.
In the original, this circle carried number one and label \emph{Stellarum Sphera Fixarum Immobilis},
meaning ``the immobile sphere of fixed stars''. This term originated in antiquity, when observes of the night sky
noticed that the lights in the sky appear to belong to two groups, those which
do not seem to move in relation to each other (\emph{stellae fixae}, fixed stars) and those
which seem to wander around (\emph{stellae errantes}, wandering stars).
 
It is worth noting that the design of the Polish coin mentioned above
had a strong competitor issued as a pattern coin in the same year. The design, carrying \emph{pr\'oba} (``pattern'') designation,
featured a more traditional portrait of Copernicus, based on the woodcut of Tobias Stimmer from before 1587.
It has not been selected for mass production, and a design with more contemporary flavour appeared on the anniversary coin.
\begin{center}
 \includegraphics[width=4.5cm]{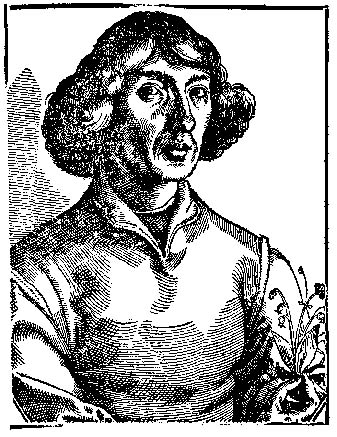} \includegraphics[width=4.5cm]{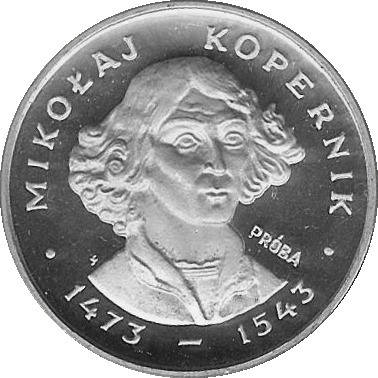} \includegraphics[width=4.5cm]{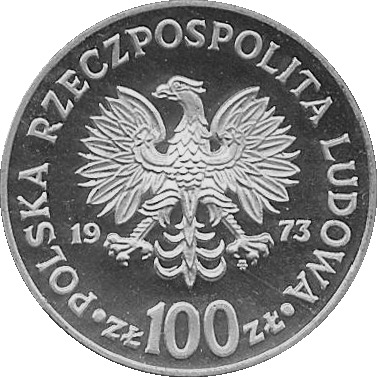}
\end{center}
In 1979, however, the same design has been used to produce .900 gold coin weighting 8 grams. Only 5000 pieces were minted,
and today it remains in very high demand.
\begin{center}
 \includegraphics[width=4.5cm]{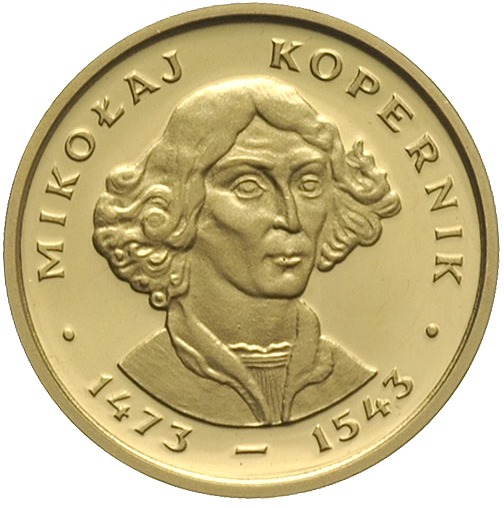} \includegraphics[width=4.5cm]{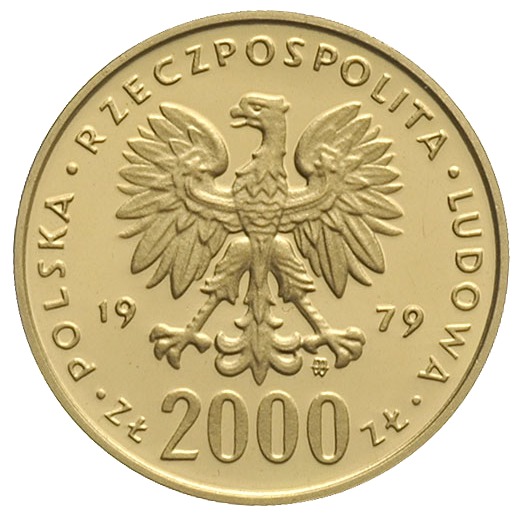}
\end{center}

The last coins with Copernicus minted before the 1994 redenomination was 5000 z\l{} coin with portrait of the astronomer
resembling anonymous 16th century painting from the Town Hall in Toru\'n. Fragment of the solar system is in the foreground,
 and panorama of Toru\'n is the background.  
\begin{center}
 \includegraphics[width=4.5cm]{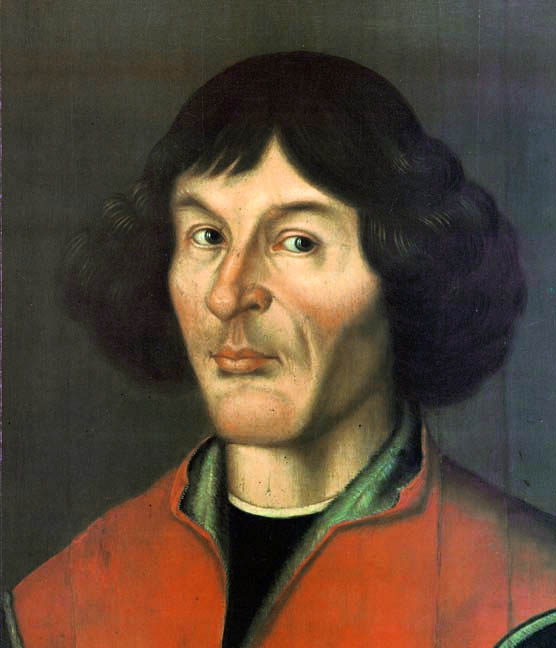} \includegraphics[width=4.5cm]{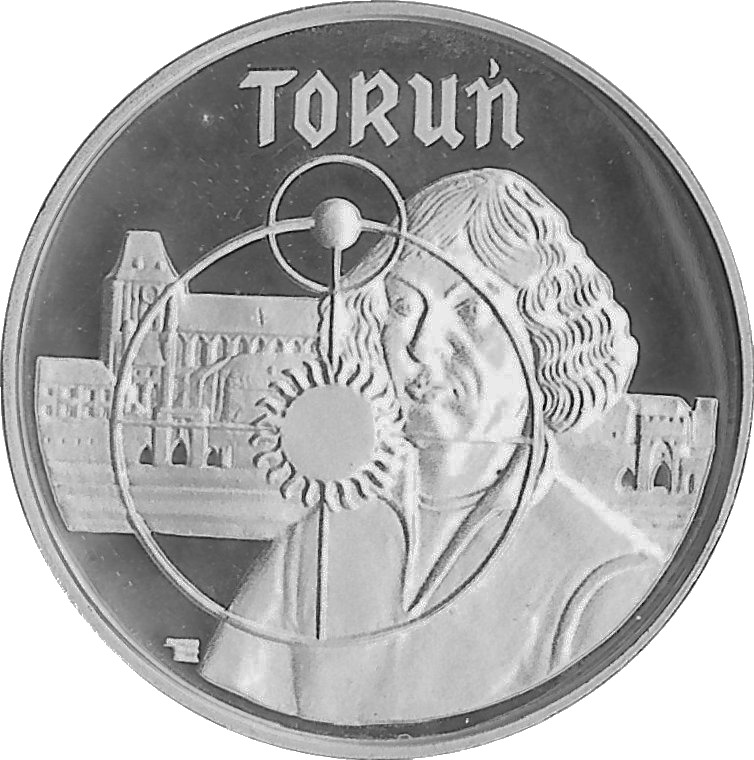} \includegraphics[width=4.5cm]{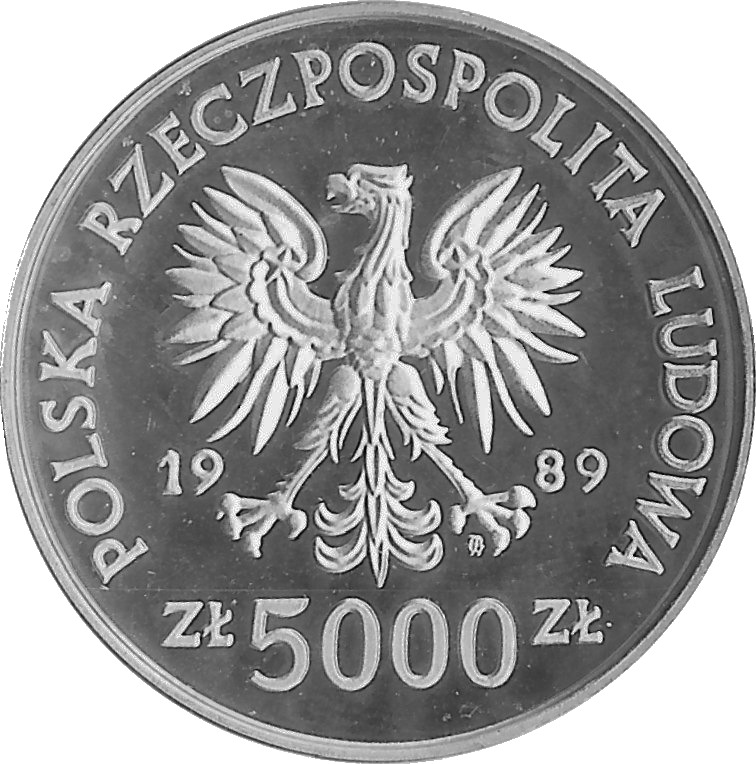}
\end{center}
The coin was made of .750 silver and had diameter of 32 mm.

Three years later, in 1992, Cook Islands issued \$10 proof coin featuring bust of Copernicus
with revolving Earth and the radiant Sun in the centre, made of one ounce of 0.925 silver.
\begin{center}
 \includegraphics[width=4.5cm]{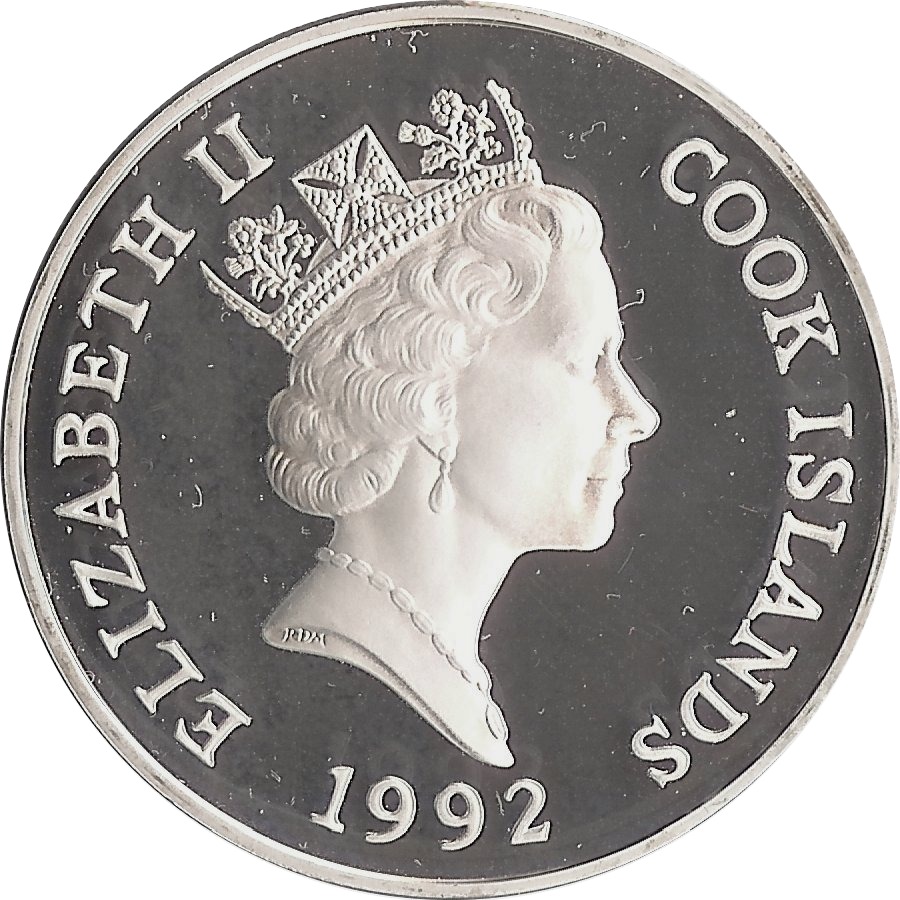} \includegraphics[width=4.5cm]{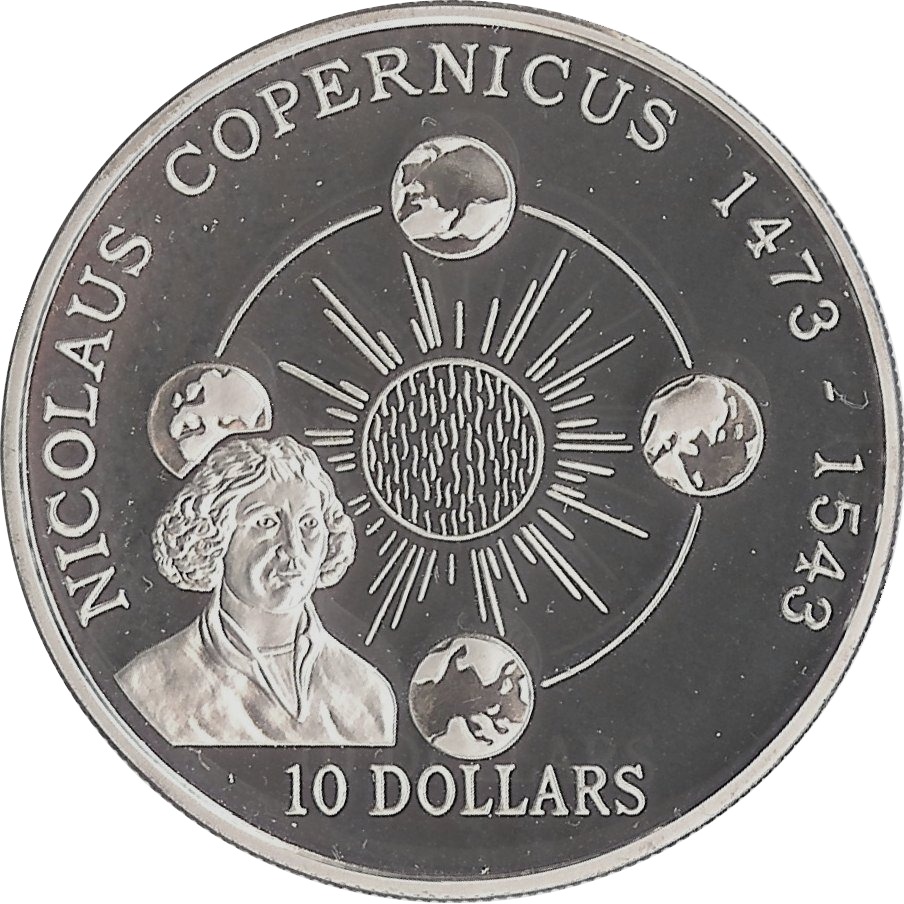}
\end{center}
Similary sized Polish coin followed in 1995, but before we describe this interesting issue, a small
digresion is in place. 

In 1989, after the fall of the communist regime, Poland again changed its official name to \emph{Rzeczpospolita Polska}
(Republic of Poland), which is evident on coins of that period. Double-digit inflation followed, and in 1995
the country had to perform  so-called redenomination, that is, introduce new \emph{z\l{}oty} equivalent to
10000 old units. The new currency is still called \emph{z\l{}oty}, but is officially denoted by PLN (the old
z\l{}oty had the ISO code PLZ).  Face values of commemorative coins once again returned to resonable
numbers, and one of the fist such commemoratives produced after the redenomination was large (38.61mm) 
 20 z\l{} coin made of one ounce of 0.925 silver. 
\begin{center}
 \includegraphics[width=4.5cm]{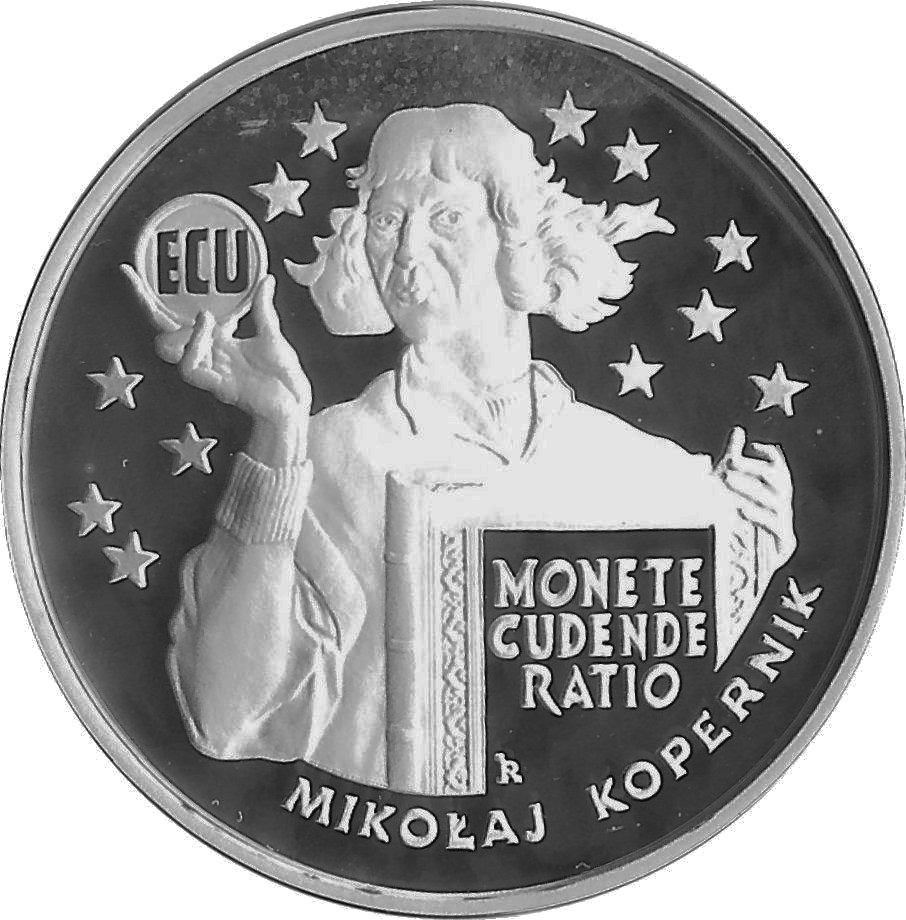} \includegraphics[width=4.5cm]{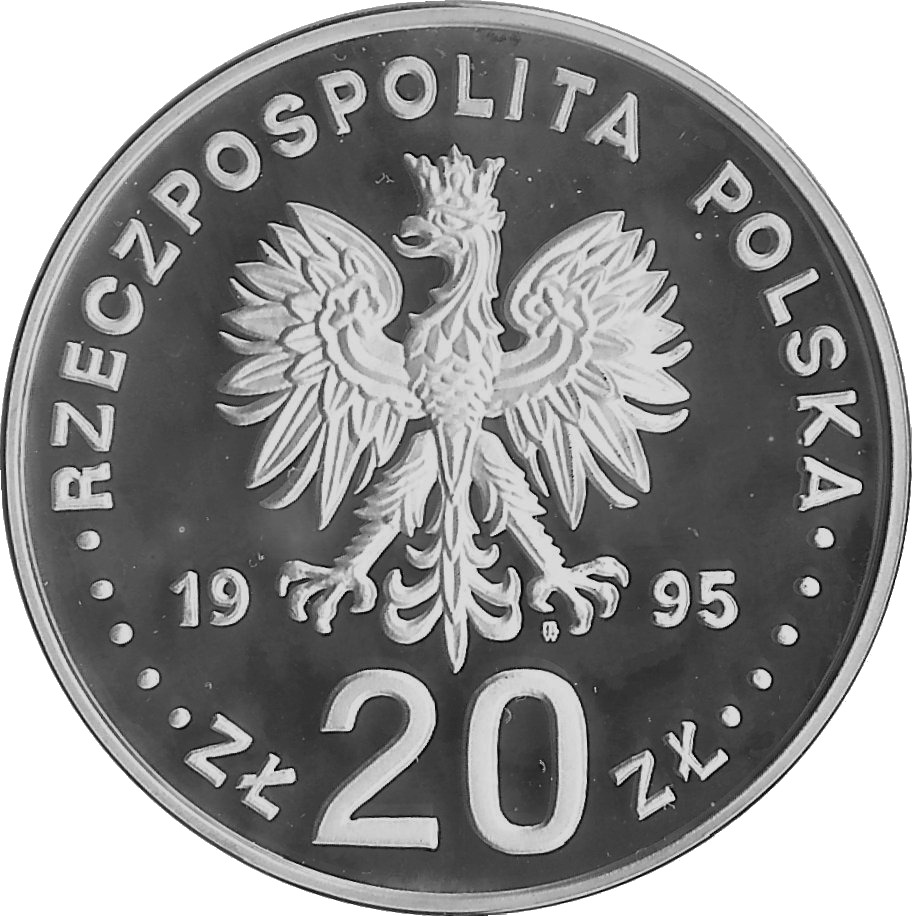}
\end{center}
Its reverse shows once again the portrait of 
Copernicus, but this time he is not pictured as an astronomer, but as an economist. In the right hand he
holds what looks like a large coin with inscription ``ECU''. This was the time when enthusiasm
for Poland's prospective  joining the European Community (now European Union)
 was running very high, and the ECU (European Currency Unit)
 was still used as the unit of account of the European Community. Twelve stars surrounding the 
head of Copernicus are also a clear reference to the flag of Europe. On the left hand side of his
silhouette we see a book titled MONETE CUDENDE RATIO. This is older (phonetic) spelling
of the title of the 1526 treaty ``On the Minting of Coin'', \emph{Monetae cudendae ratio},
a short work dealing with minting and
the quantity theory of money. Here again we must make a historical digresion.

In 1466 after the Thirteen Years' War of Poland with Teutonic Order, Poland 
recovered, after 158 years, the territory known as Royal Prussia. Initially,
the Royal Prussia was allowed to keep separate treasury and mint separate coin,
and was given rather far-reaching autonomy and privileges. It was governed by a council, 
subordinate to the Polish king, called Prussian Diet, headed by the bishop of Warmia.

At the beginning of the 16-th century, Royal Prussia was in the midst of a monetary crisis. 
Prussian coin has been suffering systematic debasing for decades, and many different coin types
were in circulation, both domestic and foreign. Most importantly,
Prussian coinage was based on a different standard that the coinage of the Polish Crown,
and for this reason, King Sigismundus the Old asked Copernicus to prepare 
a plan for monetary reform.  Copernicus complied, and in 1522 he presented his 
plan to the Prussian Diet. He proposed to uniformize  coin standards by
making  the Prussian schilling (Lat. solidus, Pol. szel{\k{a}}g) equivalent to 3 Polish units called \emph{grosze} (lat. grossus).

The plan had been accepted, although it took a couple of years to hammer out all details. In
1528 the King signed the final version of the reform bill and the production of new coins started
in the Toru\'n mint. It has been decided that a single coin would be introduced, with the face
value of 3 grosze, which later became known as \emph{trojak} (pol. for ``triple''), one of
the most characteristic and popular Polish coins, minted in various versions well until 19-th century.
It is worth noting that the Toru\'n mint produced both trojak and szal\k{a}g, although they had
identical value. We show here the solidus from 1529 and trojak from 1530, both carrying on the
obverse the inscription Sigi[smundus] Rex Polo[niae] Do[minus] Toci Pruss[iae] -- Sigismundus 
King of Poland Lord of All Prussia.

\begin{center}
 \includegraphics[width=3.5cm]{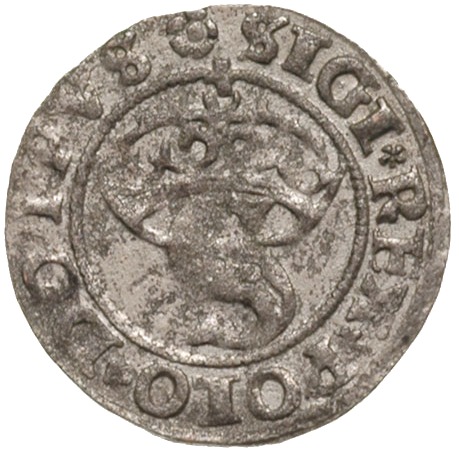} \includegraphics[width=3.5cm]{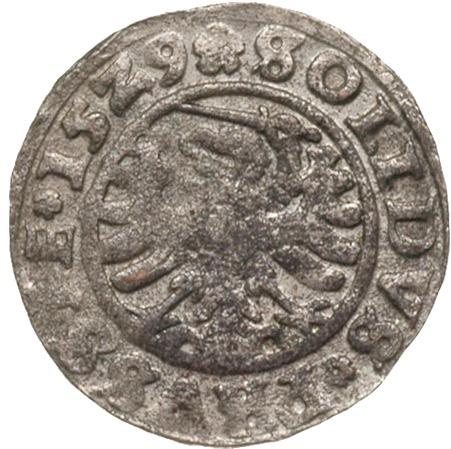}
 \includegraphics[width=3.5cm]{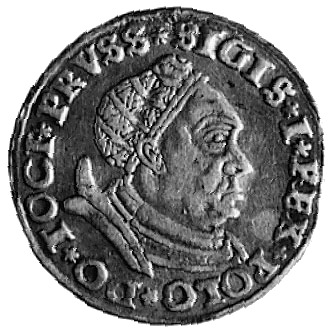} \includegraphics[width=3.5cm]{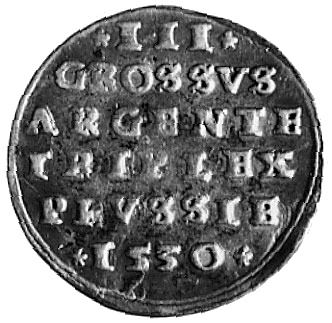}
\end{center}

In 1926, Copernicus wrote an extended version of the speech he gave to the Prussian
Diet, and gave it the title \emph{Monetae cudendae ratio}. While discussing the
nature of the monetary crisis in Prussia, he postulated in this treaty  the principle that 
``bad money drives out good''.
The principle means that when a government   undervalues one type of money,
and enforces fixed exchange rate with another type of money, 
than the ``good'' (undervalued) money will  disappear from circulation, either because it will 
leave the country, or because of hoarding. This law has been later rediscovered 
by Sir Thomas Gresham (1519--1579) and stated in his letter to Queen Elizabeth I in 1558.
As it often happens, it is now known not by the name of the original discoverer,
but as the Gresham's law. In Poland, however, it is know as Copernicus' law or sometimes
as Copernicus-Gresham's law.  It is worth stressing that the terms ``bad money'' and
``good money'' are not to be taken too literally, and it would even
be more proper to call them, respectively, ``overvalued money'' and ``undervalued money''.

\emph{Monetae cudendae ratio} makes for a very interesting reading, and its opening sentence is worth quoting
in full length:

\emph{Quanquam innumere pestes sunt quibus regna principatus, et respublice decrescere solent, haec tamen 
quatuor (meo judicio) potissime sunt:
discordia, mortalitas, terre sterilitas et monete vilitas. Tria prima adeo evidentia sunt, ut nemo ita esse
 nesciat, 
sed quartum quod ad monetam attinet a paucis et nonnisi cordatissimis consideratur, quia non uno impetu simul, 
sed paulatim, occulta quadam ratione respublicas evertit.}

\noindent This can be translated as follows:\\
``However innumerable are the scourges that used to  lead to the decline of kingdoms, 
principalities and republics, the four following are, in my judgment,  
 the most formidable: discord, pestilence, barrenness of the land, and the deterioration of the money. 
The first  three are so evident that no one is ignorant of them. But as to the fourth,
which pertains to money, with the exception of a few most prudent men, very 
few concern themselves about it; this is because it does not ruin the state at a single blow, 
but little by little, by a sort of hidden action.''

Even though this sentence was written almost 500 years ago, in the light of the financial woes
 troubling the world today, it remains  as actual as ever.

\section{Isaac Newton}
Sir Isaac Newton is quite justly considered one the foremost mathematicians and physicists 
in the history of science.  His monumental monograph \emph{Philosophiae Naturalis Principia Mathematica},
 published in 1687, established foundation of the branch of physics known today as classical
mechanics, and his discovery of the law of universal gravitation laid the groundwork for
development of modern astronomy. Newton is also a co-discoverer of  differential and integral calculus,
and author of many important discoveries in optics.

Moreover, Newton occupies a prominent place in the history of British monetary policy. 
In 1696 he obtained the post of warden of the Royal Mint, and in 1699 he became the Master of the Mint.
Although these appointments were considered to be not much more than just titles and easy salaries, Newton
actually took serious interest in his job, and became engaged in eliminating counterfeiting as well as in the coinage reform.
Here is where we encounter an interesting connection with Copernicus' \emph{Monetae cudendae ratio}.

In 1662, the old method of producing coins by hammering was superseded by a mechanical press. The coins produced
by the press were of course of much better quality that the old ones, in both shape and weight.
One could, therefore, naively expect that they would be preferred by the public and eventually replace the old hammered coins.
Yet the Copernicus-Gresham's law says otherwise. Indeed, the new coins were used for savings, export, and melting for bullion, while old,
typically underweight and worn out coins remained in circulation. The problem was compounded
by filing and clipping to obtain silver, and by 1690's it the became so severe, that a remedy had to
be found. Several solution had been proposed. The Secretary of the Treasury Lowndes recommended a recoinage at a debased standard, 
while John Locke advocated no recoinage and no debasement, proposing that the clipped coins circulate by weight,
not by the face value. Parliament adopted yet another solution, namely recoinage with no debasement, meaning that
the clipped money was to be recoined at public expense (missing silver had to be added). This operation,
known as the Great Recoinage of 1696, required massive manufacturing effort by the mint, and its 
total cost turned out to be huge, almost comparable with the total yearly tax revenue of the country.

Even though Newton was initially against the recoinage, once the
 Parliament made the decision, 
he threw all his energy into the recoinage effort. He improved efficiency of the mint, improved machinery, and 
set up branch mints. As a result, by 1698, the recoinage was practically over. 

In the following years, after becoming the the Master of the Mint, Newton advised the Treasury on monetary maters.
Among other things, he recommended lowering the price of guinea, because he observed that the price of gold
expressed in silver was higher in England than elsewhere in Europe. In December 1717, a royal proclamation has been issued
fixing the price of guinea at 21 shillings. This legislation paved the way to adopting the golds standard
a century later, and for this reason Newton is sometimes credited as the precursor of the gold standard. 
This is not quite true, as he always remained convinced that silver should be ``the true and only only 
monetary standard'' \cite{Li63}.

We will remark at this point that since the gold to silver ratio remained legally fixed in Britain after 1717, it should be
possible to use English coinage as a testbed for the Copernicus-Gresham's law, treating
gold and silver coins as two types of competing money.
 Indeed, in 2011 Thomas L. Hogan
performed detailed mathematical  analysis of historical data extracted from available financial records and concluded that
the Copernicus-Gresham's law  is ``both logically sound and empirically valid'' \cite{Hogan2011}. 

After this necessarily brief description of Newton's activities at the Royal Mint, let us
now turn to coins featuring him.
The oldest such coin  is not actually a coin in the strict sense of the word,
but rather a type of token, known as Conder tokens, named after their first cataloger, James Conder of Ipswich.
In 18-th century England, for a variety of reasons not to be discussed here, the lack of small coins was a 
perennial problem \cite{Sargent2002}. In 1787, the Parys Mining Company, who mined copper ore, saw an opportunity and 
started producing penny and halfpenny tokens, which quickly became an instant hit 
with the buying public and merchants. In years 1787--1797 almost the only small ``coins'' in circulation in 
Britain were Conder tokens. Shown below is a halfpenny Conder token from 1793 featuring somewhat unusual
profile of Sir Isaac Newton, somewhat resembling the marble bust sculpture of him
 by Louis Francois Roubiliac, now in the Wren Library in Cambridge.
\begin{center}
\includegraphics[width=4.5cm]{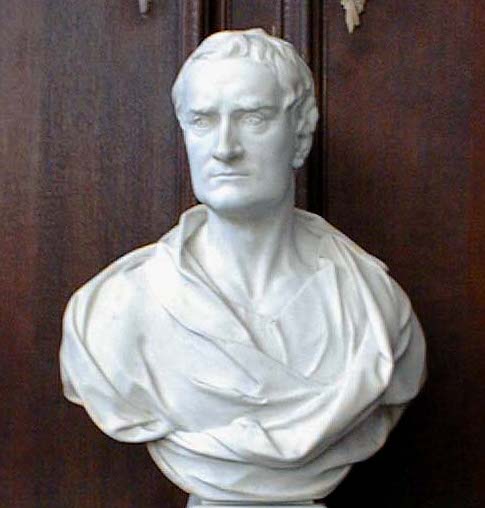}
 \includegraphics[width=4.5cm]{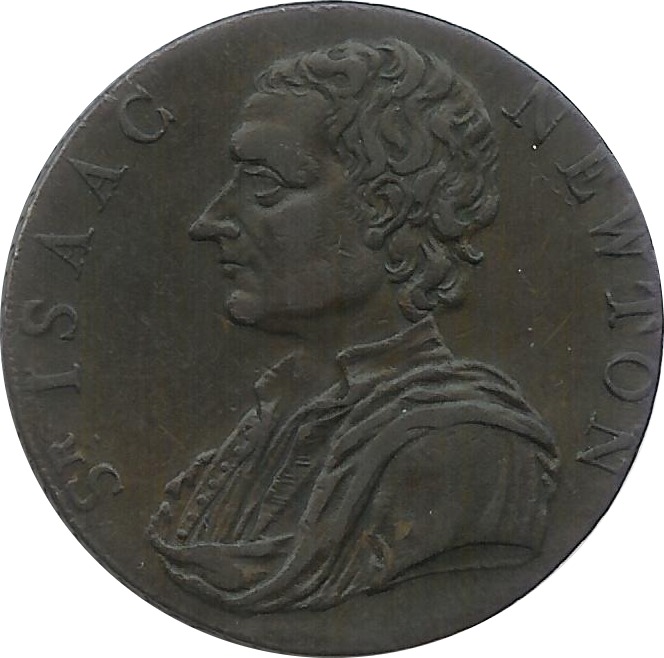} \includegraphics[width=4.5cm]{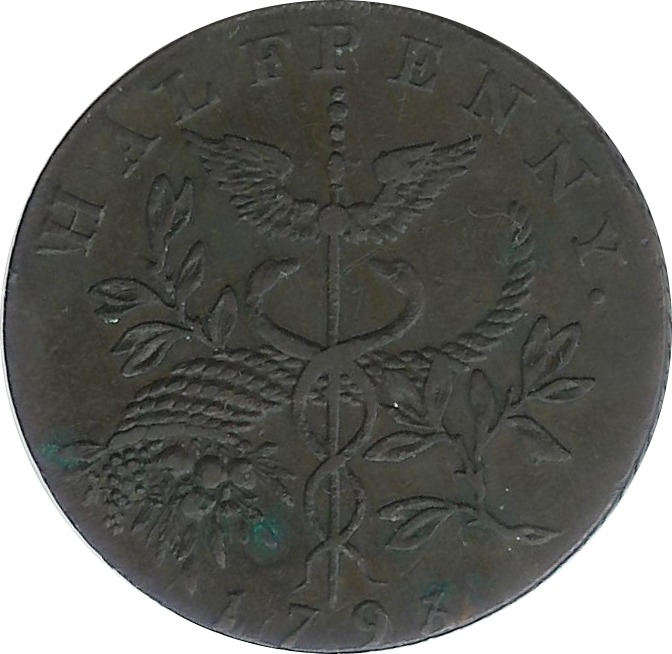}
\end{center}

Except of the aforementioned token, I am not aware of any other coin with Newton until 1993, when
Tuvalu issued 0.925 silver proof coin with denomination of \$20. On reverse we see the portrait of
Newton, and on the left schematic diagram of the planetary system with a telescope pointing toward it.  
In 1668 Newton build the first reflecting telescope, that is, a telescope which uses
a concave mirror as a primary optical element. Four years later he built an
improved version, and gave a duplicate  to the Royal Society in London. This telescope,
shown below, in now in the vaults of the the Royal Society, 
and the engraving in the on the Tuvalu represents its appearance quite well.
\begin{center}
\includegraphics[width=4.5cm]{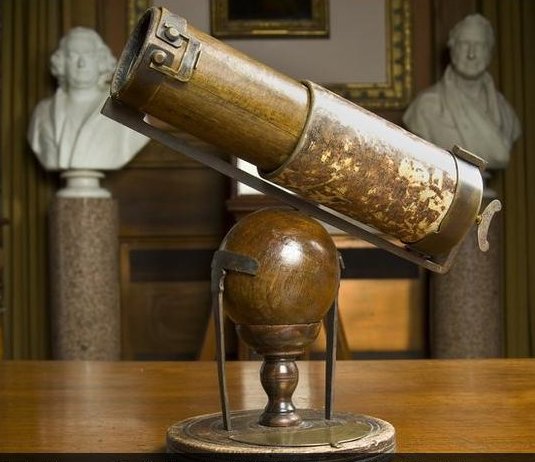} 
 \includegraphics[width=4.5cm]{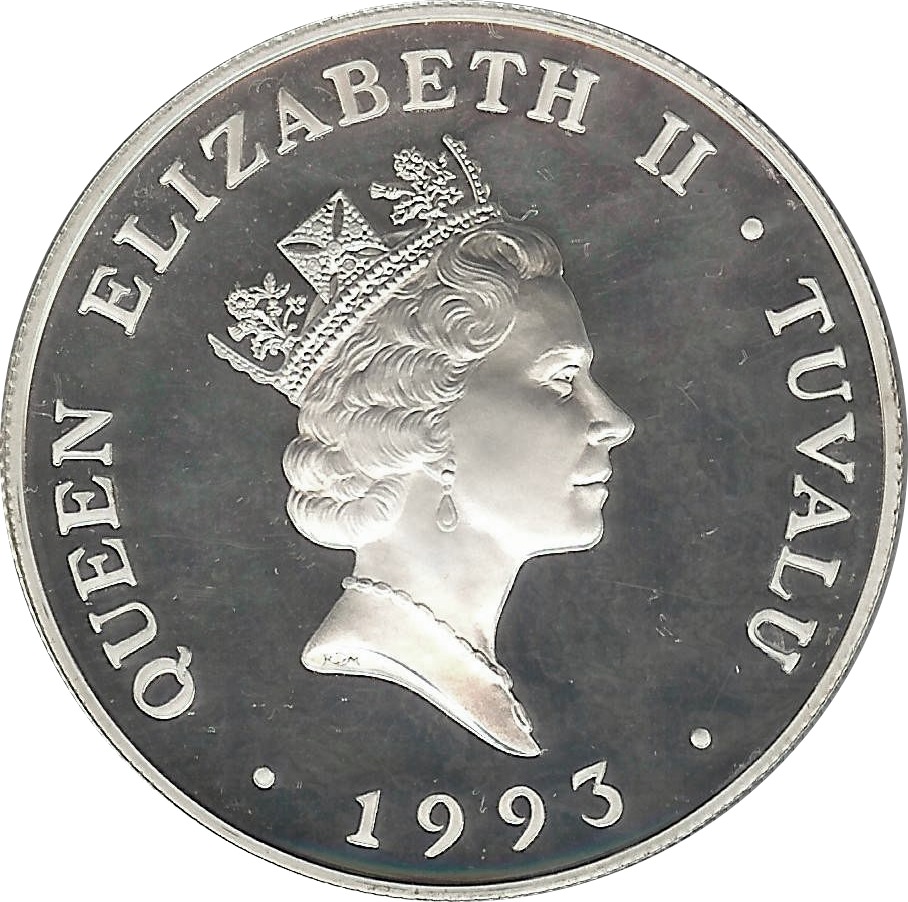} \includegraphics[width=4.5cm]{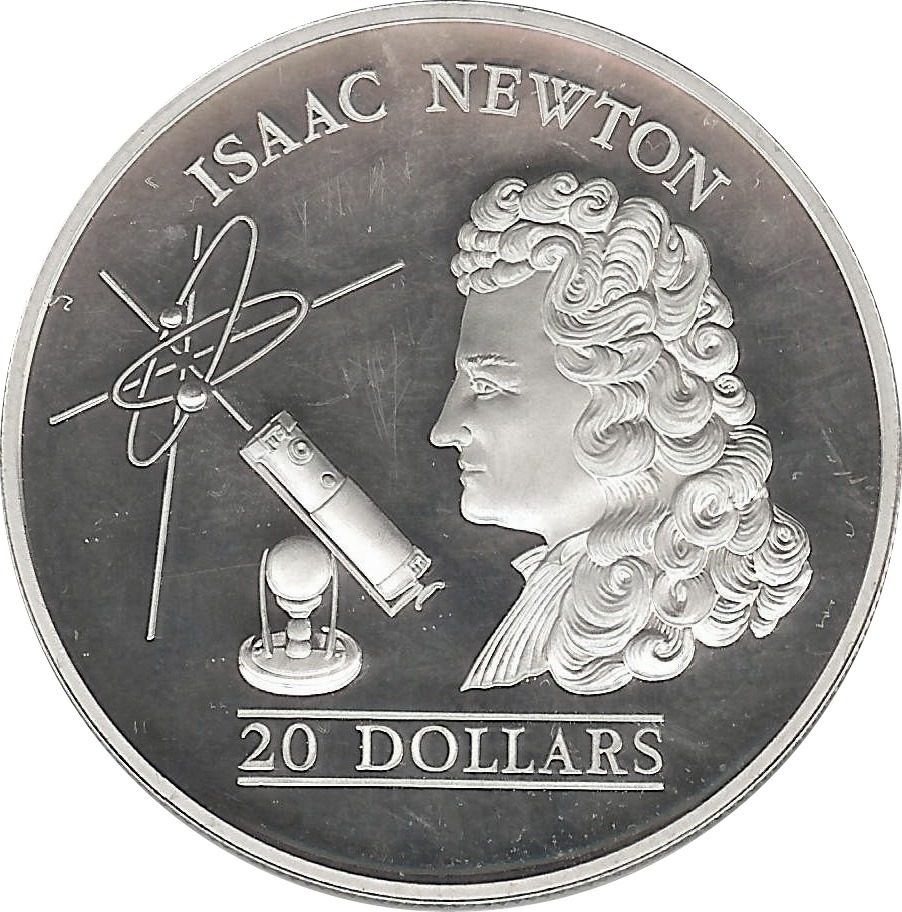}
\end{center}

Another silver proof coin comes from Cook Islands from 1997. This \$50 issue depicts Sir
Issac Newton on the reverse, and in front of him is the manuscript of \emph{Principia} as well as
the same Newtonian telescope as described above. We also see the famous apple
 which, by falling from a tree,  supposedly inspired Newton to formulate his theory of universal gravitation.
The inscription on the top reads ``17th century age of modern science -- Sir Isaac Newton's Principia''. 
\begin{center}
 \includegraphics[width=4.5cm]{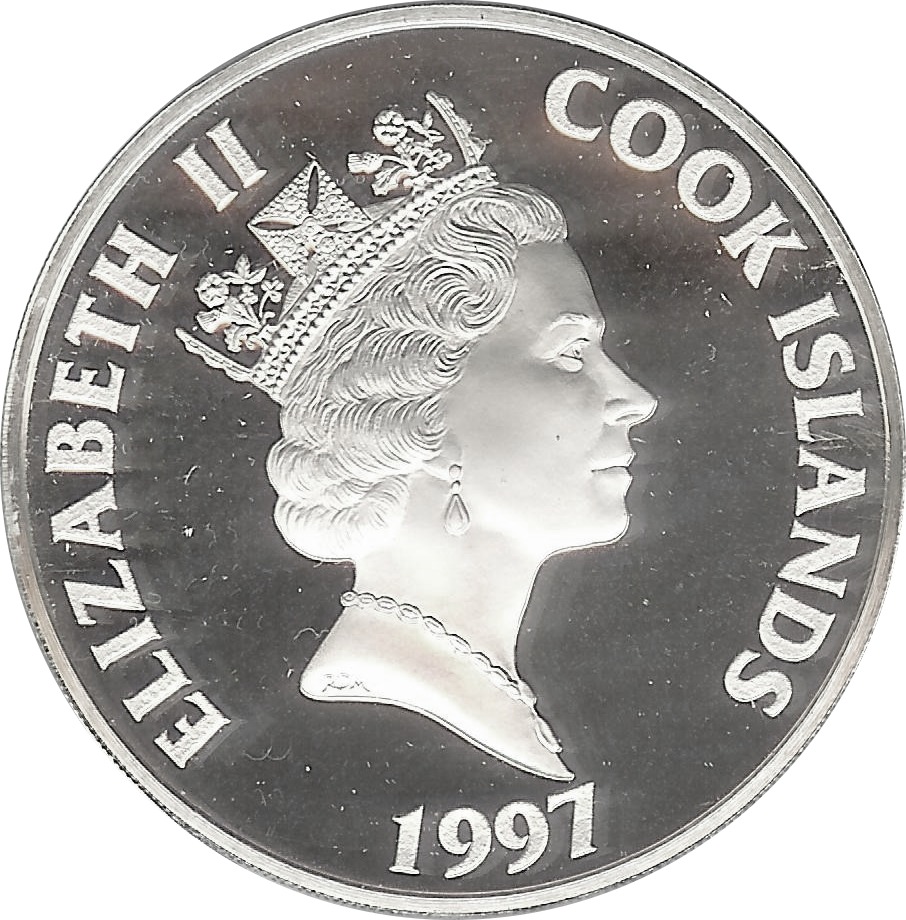} \includegraphics[width=4.5cm]{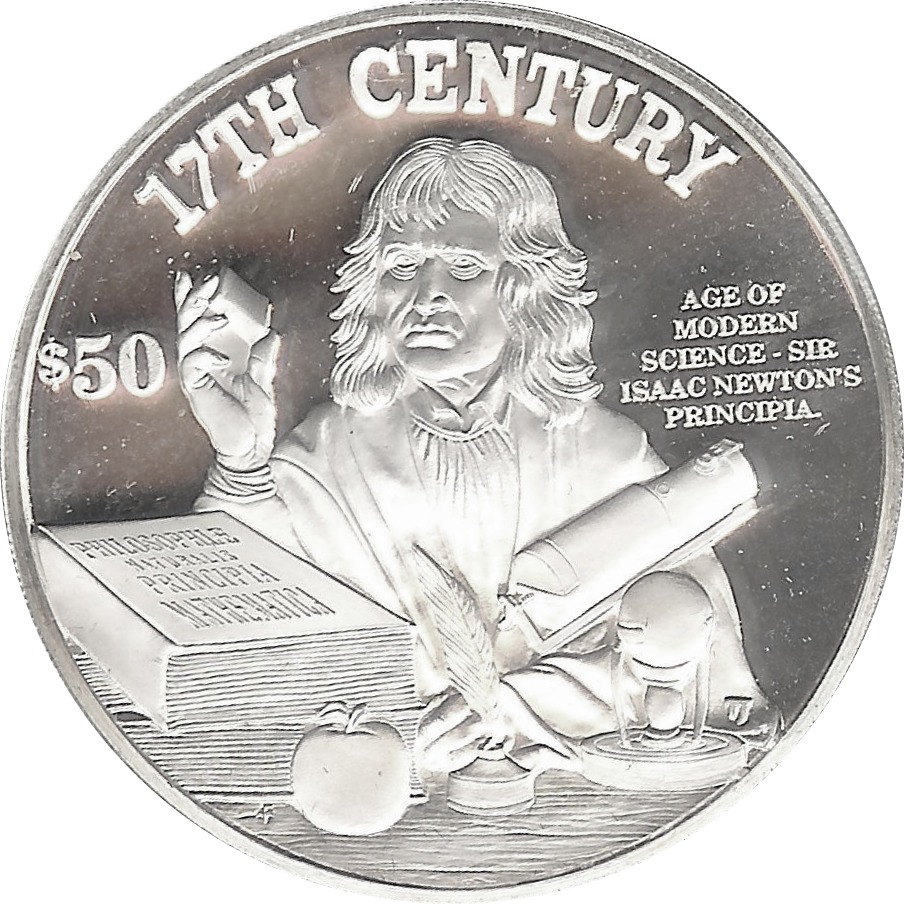}
\end{center}
The coin of Somalia from 1999 bears on the reverse the portrait of Newton  with the Newtonian
telescope superimposed below. In the middle we see the formula $F=ma=m\frac{v^2}{r}$. This is
a modern form of the second law of motion, which Newton formulated as 
\emph{Mutationem motus proportionalem esse vi motrici impressae, \& fieri secundum lineam rectam qua vis illa imprimitur.}
Today we would state it in a more expanded fashion, ``the acceleration $a$ of a body is parallel and directly 
proportional to the net force $F$ and inversely proportional to the mass $m$.'' This would
be written as $a=\frac{F}{m}$, slightly differently than what we see on the Somalian coin. The second
part of the formula on the coin is a special case of the second law, applicable to the
case when the body is moving with linear velocity $v$ around a circle of radius $r$.
\begin{center}
 \includegraphics[width=4.5cm]{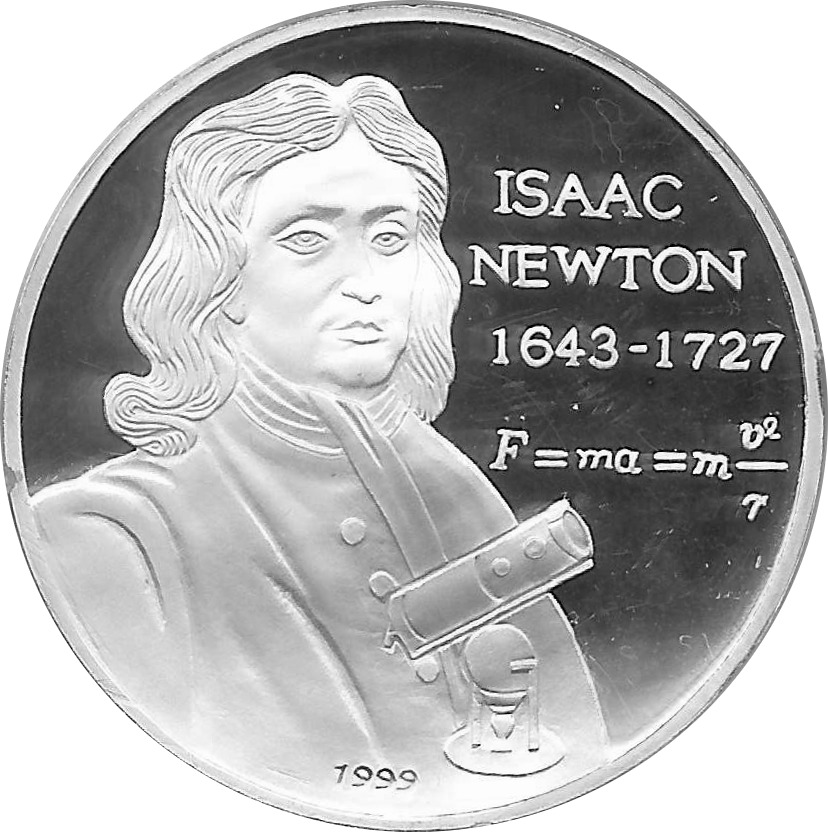} \includegraphics[width=4.5cm]{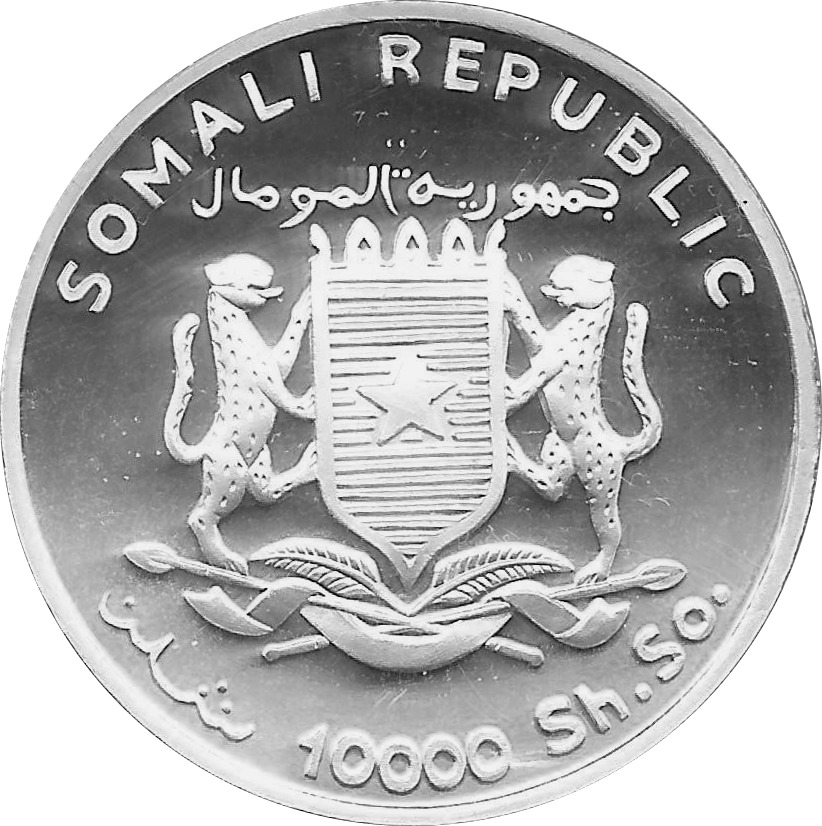}
\end{center}

The next ``Newtonian''  coin features neither his portrait nor the telescope. The reverse of the
500 Tugrik coin of Mongolia, bimetallic (silver with gold centre) issue from 1999, depicts a scheme of the solar system with
the Sun and planets from Mercury to Saturn (only these were known at the time of Newton).
Orbits of two comets are shown as well. Using Newtonian mechanics, one can show that the only possible shapes
of orbits in the so-called two body system (e.g., Sun and a planet or comet) are elliptical, parabolic, and hyperbolic.
Orbits of comets shown on the coin are thus either arcs of elongated ellipses (short-period comets), or belong to
single-apparition comets, in which case they have parabolic or hyperbolic shapes.
\begin{center}
 \includegraphics[width=4.5cm]{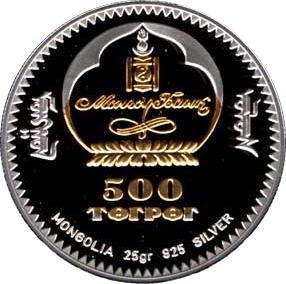} \includegraphics[width=4.5cm]{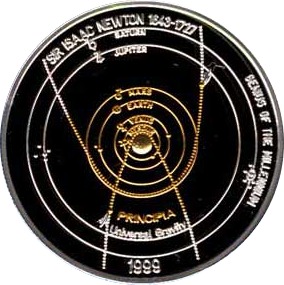}
\end{center}

We are now ready to move to the 21-st century. In 2006, the island of Alderney, a British Crown dependency,
issued an interesting coin featuring Newton. His portrait on the reverse looks very familiar, just as
 the well known 1702 portrait by Godfrey Kneller. 
\begin{center}
 \includegraphics[width=4.5cm]{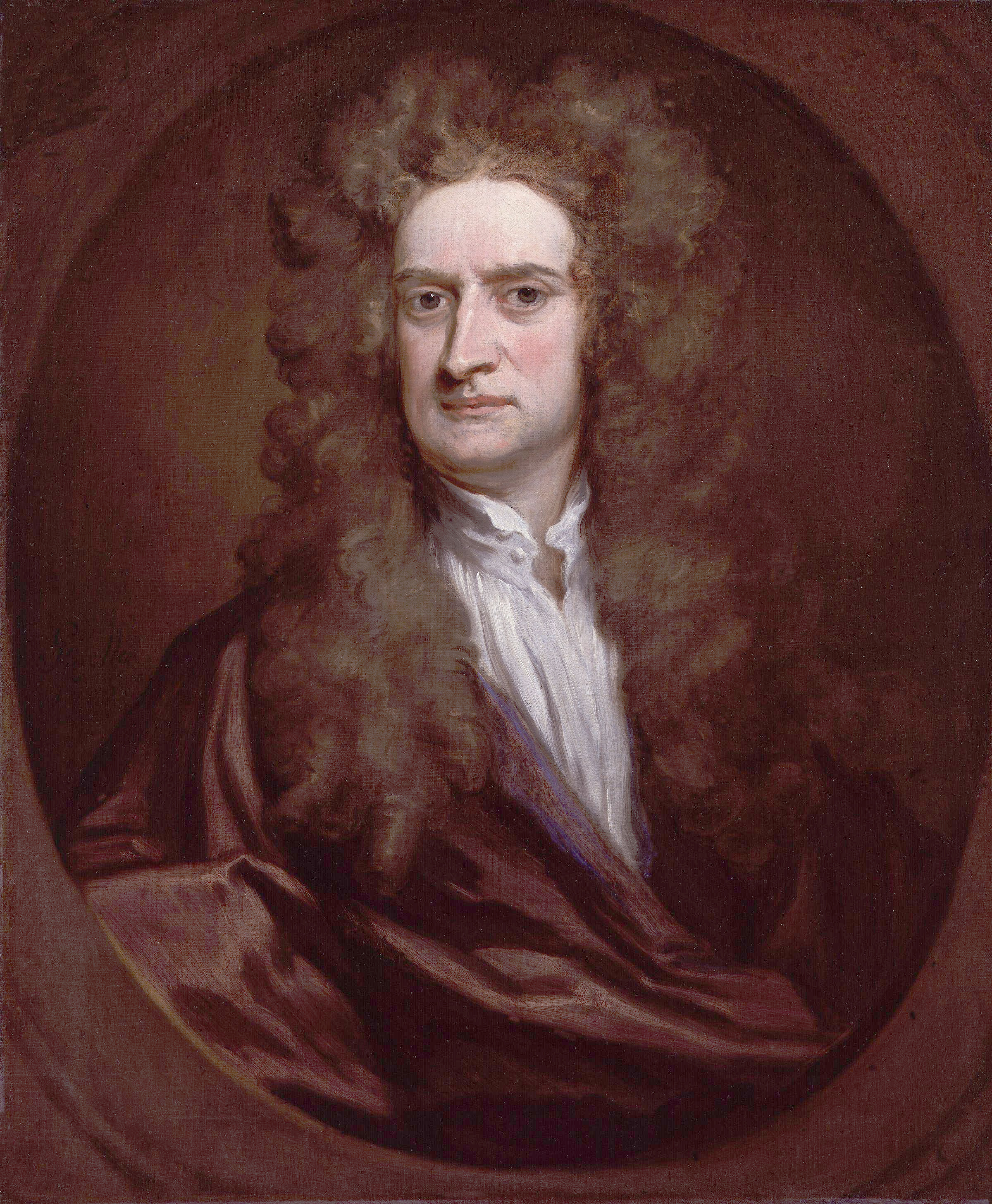}
 \includegraphics[width=4.5cm]{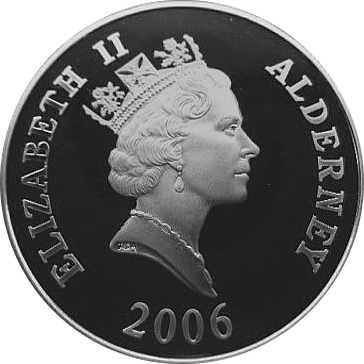} \includegraphics[width=4.5cm]{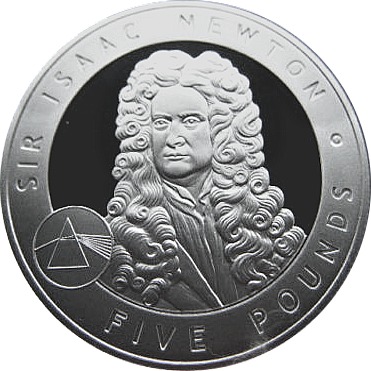}
\end{center}
The inset below shows a  triangular prism dispersing light. Newton performed this experiment in 1671, when
he introduced the word  \emph{spectrum} (Latin for ``apparition'') to describe the 
rainbow of colors emerging from the prism. His second most influential book after \emph{Principia}, entitled
\emph{Opticks} and published in 1704, discusses this experiment in detail.

After 367 years since Newton's birth, in 2009, the wait for the first official British coin commemorating
the great scientist had been finally over. As a part of the series ``Celebration of Britain'', issued by the Royal
Mint in anticipation of the London 2012 Olympic Games, 5 pound coin featuring Newton has been
designed by Shane Greeves and minted in sterling silver. The reverse displays green logo of the Games
and the inscription, ``Make not your thoughts your prison,'' taken from Shakespeare's \emph{Antony and Cleopatra}.
According to the description of the coin given by the Mint, the inscription
 ``can be seen as a message to Olympians and non-Olympians alike to strive to achieve beyond their expectations''.

The  image of Sir Isaac Newton in the background comes from the sculpture of Newton by   Eduardo Paolozzi, located in the 
piazza of the British Library. Paolozzi sculpture was inspired by the 1795 painting of William Blake, shown below.
Interestingly enough, Blake was not a fan of the Newtonian view of the universe, and his painting was not intended to be
flattering for Newton. I am not sure, therefore, what would Newton think of
 this coin.
\begin{center}
\includegraphics[width=5.5cm]{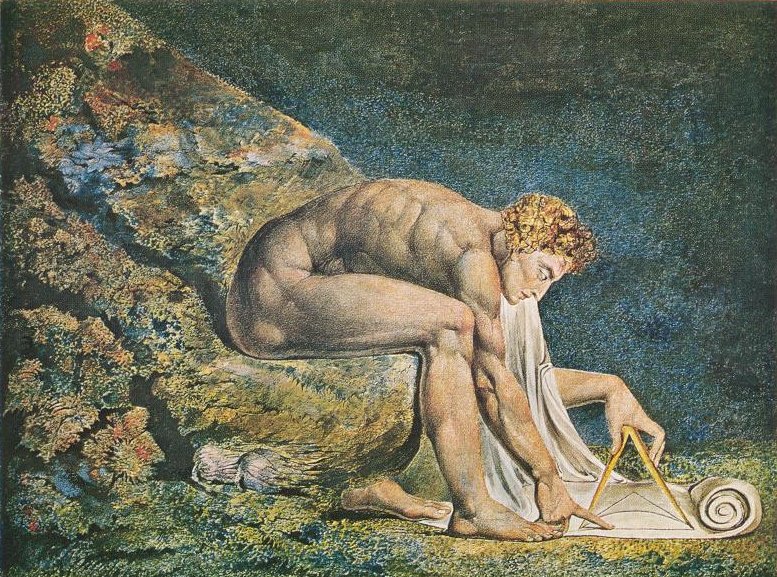} \includegraphics[width=5.5cm]{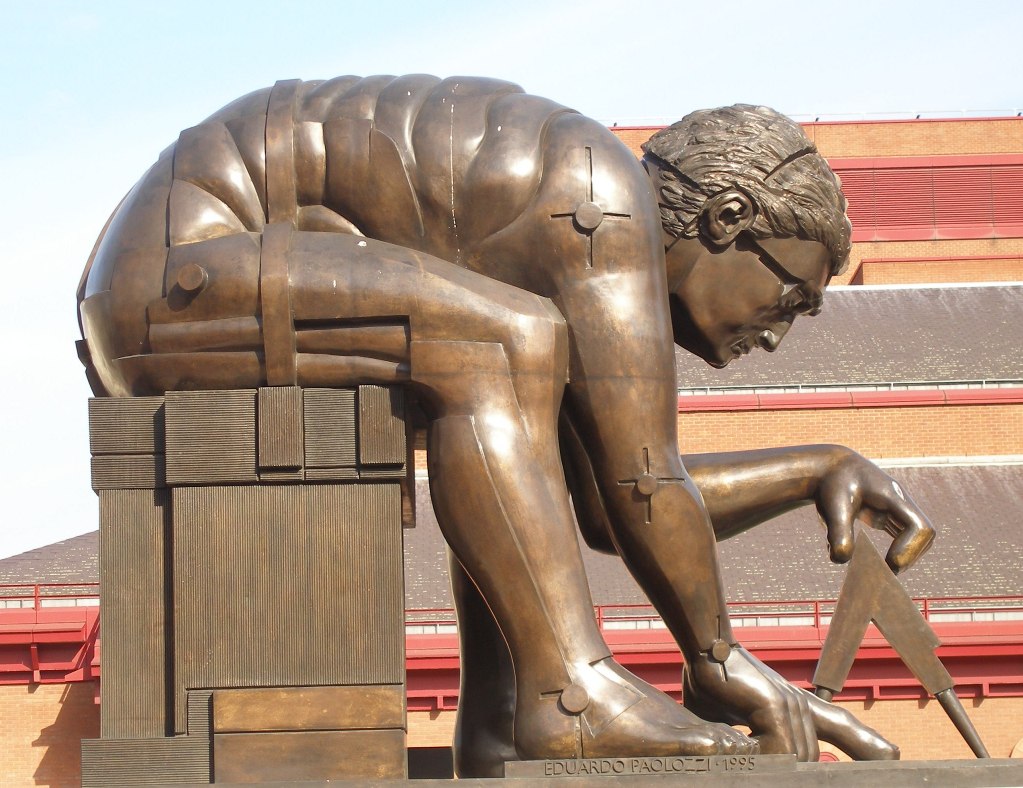} \\
 \includegraphics[width=4.55cm]{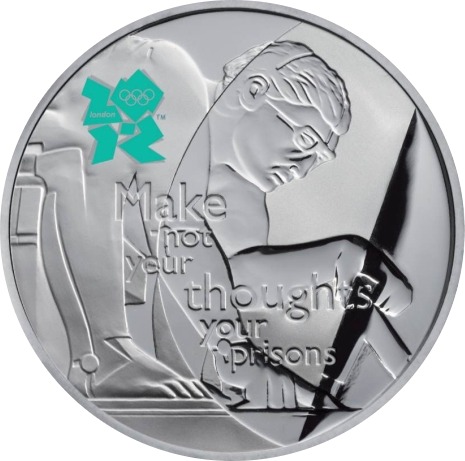} \includegraphics[width=4.5cm]{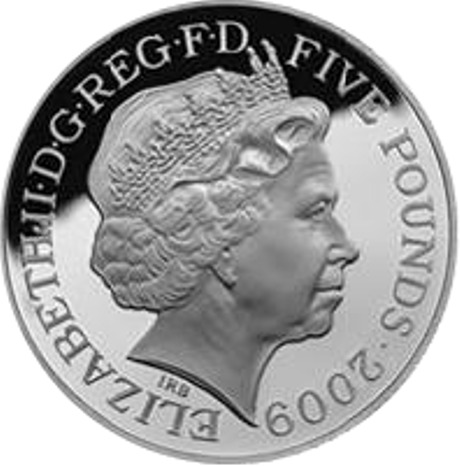}
\end{center}
But regardless of what he would think, one can make a sad observation.
In  today's Britain, similarly as in Canada, sports seem to be valued much more than intellectual activities, in particular 
more than
mathematics. The Royal Mint finally commemorated Newton, but only in passing, as a background for promotion
of Olympic Games. The policy of the Royal Canadian Mint is, unfortunately, no different: alongside of countless coins featuring
Olympic Games,  hockey players,  and NHL teams, you will  not find a single Canadian coin even remotely related to 
mathematics or mathematical sciences.\footnote{Post-publication update: in 2015, three years after the publication of this article, RCM issued 10oz silver
coin commemorating Albert Einstein. Also in 2015, silver and gold bullion coins with $E=mc^2$ privy mark were issued.}
\vskip 1em
{\small 
\noindent \textbf{Photo credits:} author's collection,
wbcc-online.com, and royalmint.com.
 Solidus, trojak, and 1925  pattern from
the archive of Warsaw Numismatic Centre (www.wnc.pl), reproduced with permission.} 

\nocite{Li63,Sargent2002,Anthony1960,Kalkowski1981}

\providecommand{\href}[2]{#2}\begingroup\raggedright\endgroup
\end{document}